%% file: m4-9.tex
\documentclass{gtart}

\input gtmonout

\volumenumber{4}
\volumename{Invariants of knots and 3-manifolds (Kyoto 2001)}
\volumeyear{2002}
\papernumber{9}
\pagenumbers{119}{141}
\received{16 December 2001}
\revised{14 April 2002}
\accepted{22 July 2002}
\published{19 September 2002}

\usepackage{amsmath,amssymb, epsf}    

\def\Tors{{\rm Tors}}
\def\btm{\begin{thm}}
\def\etm{\end{thm}}
\def\Spl{{\rm Sp}}
\def\so{{\mathfrak s}{\mathfrak o}}
\newcommand{\supp}[1]{{\scriptscriptstyle [\leq #1]}}
\def \lz  {\langle}
\def \rz  {\rangle}
\def \fr{ \mbox{{\it\large f}}}
\def\sel{{\mathfrak s}{\mathfrak l}}
\def\sp{{\mathfrak s}{\mathfrak p}}
\let\ep\endproof

\newcommand{\vs}[1]{}

\newtheorem{theorem}{Theorem}
\newtheorem{corollary}{Corollary}

\newtheorem{conjecture}[theorem]{Conjecture}

\newtheorem{lemma}[theorem]{Lemma}
\theoremstyle{definition}
\newtheorem{definition}{Definition}

\newtheorem{example}{Example}

\def\dashfill{\cleaders\hbox{$\mkern2.2mu \_ \mkern2.2mu$}\hfill}

\def\underdash#1{\mathop
{\vtop{\ialign{##\crcr
$\hfil\displaystyle{#1}\hfil$\crcr\noalign{\kern.5pt\nointerlineskip}
\dashfill\crcr\noalign{\kern.5pt}}}}\limits}

\def\rightarr{\mbox{\begin{picture}(15,1)(0,0)\put(0,0)
{\vector(1,0){15}}\end{picture}}}

\def\underrightarrow#1{\mathop
{\vtop{\ialign{##\crcr
$\hfil\displaystyle{#1}\hfil$\crcr\noalign{\kern1pt\nointerlineskip}
\rightarr\crcr\noalign{\kern.5pt}}}}\limits}

\def\leftarr{\mbox{\begin{picture}(15,1)(0,0)\put(15,0)
{\vector(-1,0){15}}\end{picture}}}

\def\underleftarrow#1{\mathop
{\vtop{\ialign{##\crcr
$\hfil\displaystyle{#1}\hfil$\crcr\noalign{\kern1pt\nointerlineskip}
\leftarr\crcr\noalign{\kern.5pt}}}}\limits}

\def\D{\displaystyle}

\newcommand{\ext}[1] {\mbox{\raisebox{.4ex}{$\bigwedge^{\!#1}$}}\mkern-1mu}

\newcommand{\Z}{{\mathbb Z}}
\newcommand{\R}{{\mathbb R}}
\newcommand{\N}{{\mathbb N}}

\newcommand{\F}{{\mathbb F}}

\newcommand{\M}{{\mathbb M}}

\newcommand{\posit}[2]
{\raise -1.4ex\hbox{${\textstyle #1}\atop {\stackrel{\uparrow}{#2}}$}}

\def \bm {\boldmath}
\def \ubm {\unboldmath}

\newcommand{\bmu}
{\mbox{\bm$\mu$\ubm}}
\newcommand{\bnu}
{\mbox{\bm$\nu$\ubm}}
\newcommand{\bla}
{\mbox{\bm$\lambda$\ubm}}
\newcommand{\bka}
{\mbox{\bm$\kappa$\ubm}}
\newcommand{\brh}
{\mbox{\bm$\rho$\ubm}}
\newcommand{\bxi}
{\mbox{\bm$\xi$\ubm}}
\newcommand{\bOm}
{\mbox{\bm$\vec o$\ubm}}

\def\symm{\leftrightarrow}

\def\opc{{\scriptstyle =}}

\def\drulefill{$\mathord\opc\mkern-3mu\cleaders
                  \hbox{$\mkern-2mu\mathord\opc\mkern-2mu$}\hfill\mkern-2mu \mathord\opc$}

\def\dov#1{\vbox{\ialign{##\crcr\noalign
             {\kern-3pt\nointerlineskip}\drulefill \crcr\noalign
             {\kern0pt\nointerlineskip}
             $\hfil\displaystyle{#1}\hfil$\crcr}}}

\def\trace{\mbox{\rm\small trace}}
\def \into {\hookrightarrow}
\def\onto {-\!\!\!\!\!{\to\mkern-14mu\to}}
\def\ben{\begin{enumerate}}
\def\een{\end{enumerate}}

\begin{document}

\title{$p$-Modular TQFT's, Milnor torsion\\and the Casson-Lescop
invariant}
\asciititle{p-Modular TQFT's, Milnor torsion\\and the Casson-Lescop
invariant}
\author{Thomas Kerler}
\address{Department of Mathematics, The Ohio State University\\231 
West 18th Avenue, Columbus OH 43210, USA}
\email{kerler@math.ohio-state.edu}
\primaryclass{57R56, 57M27}
\secondaryclass{17B10, 17B37, 17B50, 20C20}
\begin{abstract}
We derive formulae which lend themselves to TQFT interpretations of
the Milnor torsion, the Lescop invariant, the Casson invariant, and
the Casson-Morita cocyle of a 3-manifold, and, furthermore, relate
them to the Reshetikhin-Turaev theory.
\end{abstract}

\keywords{Milnor-Turaev torsion, Alexander polynomial,
Casson-Walker-Lescop invariant, Casson-Morita cocycle, TQFT, 
Frohman-Nicas theory, Reshetikhin-Turaev theory, $p$-modular
representations}

\asciikeywords{Milnor-Turaev torsion, Alexander polynomial,
Casson-Walker-Lescop invariant, Casson-Morita cocycle, TQFT, 
Frohman-Nicas theory, Reshetikhin-Turaev theory, p-modular
representations}

\maketitlepage

\addtocounter{section}{-1}

\section{Introduction and summary}

Invariants of 3-manifolds that admit extensions to topological quantum 
field theories (TQFT) are structurally highly organized. Consequently, 
their evaluations permit 
an equally deeper insight into the topological structure of the underlying 3-manifolds
beyond the mere distinction of their homeomorphism types. 
Although the notion and many examples of TQFT's have been around for 
more than a decade there are still surprisingly large gaps in the understanding 
of the explicit TQFT content of  the ``classical''  
invariants such as the Milnor-Turaev torsion or the Casson-Walker-Lescop invariant.

In this article we derive formulae which lend themselves 
to TQFT interpretations of the 
Milnor torsion, the Lescop invariant, the Casson invariant,
and the Casson-Morita cocyle  of a 3-manifold. Specifically, these invariants 
are expressed in (\ref{eq-AlexPoly}) of Theorem~\ref{thm-LefTr}, in (\ref{eq-LesMon})
of Theorem~\ref{thm-LesMon}, in (\ref{eq-Cassmat}) of Theorem~\ref{thm-Cassmat},  
and in (\ref{eq-cocycmat}) of Theorem~\ref{thm-cocycmat} as traces and matrix elements
of operators acting on  $\ext *H_1(\Sigma)$ for a surface $\Sigma$. 
We relate these formulae to previous results in \cite{Ke00}, \cite{Ke01}, \cite{KeF}
and \cite{GK} on the Frohman-Nicas and Reshetikhin-Turaev theories. In the course, we develop 
the general notion of a {\em $q/l$-solvable TQFT} and consider reductions to the $p$-modular cases,
as needed for the quantum theories. As an example for the additional structural depth that
TQFT interpretations provide  we describe results from \cite{GK} that allow us to read 
the cut-numbers of  3-manifolds from the coefficients of the Alexander polynomial, using 
their relations with the Reshetikhin-Turaev invariants.

After  a review of the Alexander polynomial and the Milnor-Turaev
torsion in Section~1 we 
introduce  in Section~2 the Frohman-Nicas TQFT modeled over the
$\Z$-cohomology of the surface Jacobians, as well as 
its Lefschetz components. The latter lead us to define the fundamental
torsion weights $\Delta_{\varphi}^{(j)}(M)\in\Z\,$ for a 3-manifold, 
$M$, together with a 1-cocycle, $\varphi$, which turn out to be representation
theoretically more adequate recombinations of the coefficients of the 
Alexander polynomial. We derive an expression of the Lescop invariant $\lambda_L$
in terms of these fundamental weights. In Section~3 
we construct the  Johnson-Morita extensions ${\cal U}^{(j)}$
of pairs of Lefschetz components, which have degrees differing by 3. 
We review elements of Morita's theory for the Casson invariant
$\lambda_C$ and the Casson-Morita cocyle $\delta_C$, and 
derive formulae that express $\lambda_C$ and $\delta_C$ in terms of 
matrix elements of  operators acting via exterior multiplication on
the cohomology of the Jacobian of a surface. We further introduce in
Section~3 the notion of a 1/1-solvable TQFT over  ${\mathbb M}[{\sf
y}]/{\sf y}^2$ for a commutative ring $\mathbb M$, which has an
obvious generalization to $q$/$l$-TQFT's and can be viewed as 
a lowest order deformation theory for the ${\cal U}^{(j)}$. We derive some 
immediate functorial properties, and use these  
and the  results for $\lambda_C$ and $\delta_C$ 
to infer criteria for when a such a TQFT realizes the Casson invariant. 
In Section~4  we discuss the modular structure and characters of reductions 
of the TQFT's  ${\cal U}^{(j)}$ from Section~2 into $\F_p=\Z/p\Z$. 
Finally, in Section~5, 
we give an example of a 1/1-solvable TQFT
over $\F_5[{\sf y}]/{\sf y}^2$ using the Reshetikhin Turaev theory at 
a 5-th root of unity. As applications we
 show that for a 3-manifold with $b_1(M)\geq 1$
the quantum invariant is given in lowest order by the Lescop
invariant, and describe how to obtain from the structure theory of
these TQFT methods for computing the cut-number of a 3-manifold.  

\section{The Alexander polynomial, Reidemeister-Milnor-Turaev 
torsion and the Lescop invariant}

We start with a short review of 
the Alexander Polynomial
$\Delta_{\varphi}(M)\in \Z[\Pi]/\pm\Pi\,$, 
which is defined  for a compact, connected, oriented
3-manifold $M$ with $b_1(M)\geq 1$, and  an epimorphism 
$\varphi:H_1(M,\Z)\onto\Pi$, where $\Pi$ is a 
free abelian group of rank $k\leq b_1(M)$. 

The map $\varphi$
defines a covering space
$\widetilde M_{\varphi}\to M$, such that we have a short exact sequence 
$1\to \pi_1(\widetilde M_{\varphi})\to \pi_1(M)\to\Pi\to 1$ and 
$\widetilde M_{\varphi}$ admits a $\Pi$-action by Deck-transformations. 
In particular, $H_1(\widetilde M_{\varphi},\Z)$ admits a $\Pi$-action and
thus becomes a $\,\Z[\Pi]\,\cong\,$ $\Z[t_1^{\pm 1}, \ldots, t_k^{\pm 1}]$ 
module. The Alexander polynomial 
 is then defined as the generator of
the smallest principal ideal containing  the first elementary ideal 
of $H_1(\widetilde M_{\varphi},\Z)$. In the case where $k=b_1(M)$,
which means 
$\Pi=H_1^{\it fr}(M)=H_1(M)/\Tors(H_1(M))\,$, we simply write $\Delta(M)$
without subscript. 

Let us briefly explore this definition in the rank 
$k=1$ case, meaning  $\Pi=\Z$. Quite often the $\Z[t^{\pm 1}]$-module  admits a  
presentation of the form $H_1(\widetilde M_{\varphi},\Z)\cong \Z[t^{\pm 1}]^m/A(\Z[t^{\pm 1}]^m)$ 
obtained from an $m\times m$ Alexander Matrix $A$ with coefficients in $\Z[t^{\pm 1}]$.
The Alexander polynomial  is thus represented by 
$\Delta_{\varphi}(M)=\pm t^l\cdot \det(A)$, with $l\in\Z$. One way of computing 
such an Alexander Matrix $A$ is as follows. 
For $\varphi:H_1(M)\onto \Z$ we find a dual, two-sided, embedded surface 
$\Sigma\subset M$ of some genus $g$. 
If we remove a collar neighborhood of $\Sigma$ we
obtain $C_{\Sigma}=M-N(\Sigma)$, which we view as
a cobordism from $\Sigma$ to itself. 
Let $i^{\pm}:\Sigma\into C_{\Sigma}$ be the inclusion maps onto the upper and the lower
boundary component of $C_{\Sigma}$.
Denote by $A^{\pm\it fr}: H_1(\Sigma_g)\to H^{\it fr}_1(C_{\Sigma})$
the maps
induced by $H_1(i^{\pm})$ onto the free parts of the homology groups.
 The sign convention should be such that
$A^{\pm\it fr}= id$ if $M=S^1\times \Sigma$ and 
$C_{\Sigma}=[0,1]\times\Sigma$. The definitions imply the formula 
\begin{equation}\label{eq-AlexForm}
\Delta_{\varphi}(M)\;\;=\;\;\pm t^{-g}\cdot |\Tors(C_{\Sigma})|\cdot 
\det(A^{+\it fr}-tA^{-\it fr})\;\;.
\end{equation}
We assume in (\ref{eq-AlexForm}) that $b_1(C_{\Sigma})=2g$ so that the $A^{\pm \it fr}$
are indeed $(2g)\times (2g)$-matrices. For example, if $M$ is the mapping torus 
with gluing map
$\psi$ and $\varphi$ is canonical then $\Delta_{\varphi}(M)$ is the characteristic polynomial of
$H_1(\psi)$. The case 
$b_1(C_{\Sigma})>2g$ is equivalent to saying that $\varphi$ factors through
an epimorphism $\pi_1(M)\onto \Z*\Z\onto \Z$ onto the (non-abelian)
free group in two generators, and implies that $\Delta_{\varphi}(M)=0$.
The additional factor in (\ref{eq-AlexForm}) and Poincar\'e duality \cite{Mi61} yield the 
symmetrized version of the Alexander polynomial. Moreover, we choose the sign such that
$\Delta_{\varphi}(1)\geq 0$. A straight forward homological computation shows 
$\Delta_{\varphi}(M)(1)= |\Tors(H_1(M))|$ if $b_1(M)=1$ 
and $\Delta_{\varphi}(M)(1)=0$ if $b_1(M)\geq 2$. If ${\cal K}\subset S^3$
is a knot we obtain the usual 
Alexander polynomial $\Delta_{\cal K}$ of knot theory by applying
the above either to the knot complement or to the 3-manifold obtained by
doing 0-surgery along ${\cal K}$. 

A closely related  invariant is the Reidemeister-Milnor-Turaev torsion 
$\tau_{\varphi}(M)\in\Z[\Pi]$. It is 
defined as the Reidemeister torsion of the acyclic chain complex associated
to the local system defined by $\varphi$, see \cite{Mi61}, \cite{Tu75} for details.
Again, for $\Pi=H_1^{\it fr}(M)$ we
simply write $\tau(M)$. In Theorems
A and B, and 4.1.II. of \cite{Tu75} Turaev shows the following relations.

\begin{theorem}\label{thm-turtor}{\rm\cite{Tu75}}\qua
Let $\varphi:H_1(M)\onto \Pi$ with $k={\rm rank}(\Pi)$: 
\ben
\item If $k\geq 2$ then $\tau_{\varphi}(M)=\Delta_{\varphi}(M)$
(so if $b_1(M)\geq 2$ then $\tau(M)=\Delta(M)$). 
\item If $k=1$ and $\partial M = \emptyset$ then 
$\displaystyle \tau_{\varphi}(M)\,=\, {(t-1)}^{-2}\Delta_{\varphi}(M)$.
\item If $b_1(M)\geq 2$, $k=1$, and $\partial M = \emptyset$
then $\Delta_{\varphi}(M)=\overline{\varphi}(\Delta(M))(t-1)^2t^{-1}$. 
\item If $k=1$ and  $\partial M\neq \emptyset$ then 
$\displaystyle \tau_{\varphi}(M)\,=\, {(t-1)}^{-1}\Delta_{\varphi}(M)$.
\item If $b_1(M)\geq 2$, $k=1$, and $\partial M \neq \emptyset$
then $\Delta_{\varphi}(M)=\overline{\varphi}(\Delta(M))(t-1)$.
\een 
\end{theorem}

We next recall the relation between the Alexander polynomial of a link and 
that of the corresponding closed 3-manifold obtained by surgery. 
Let ${\cal L}\subset N$ be a framed link
in a connected rational homology sphere $N$  with $n\geq 2$ number of
components, 
all of which
have 0 framings and 0 linking numbers. Denote by $N_{\cal L}$ the 3-manifold 
obtained by surgery along $\cal L$ so that $b_1(N_{\cal L})=n$.
Hence, $\Delta(N_{\cal L})\in\Z[t_1^{\pm 1},\ldots,t_n^{\pm 1}]$, 
where the generator $t_j$ is given by the meridian $l_j$ of the 
$j$-th component of ${\cal L}$. By Theorem E of \cite{Tu76} and (1) of 
Theorem~\ref{thm-turtor} we have $\Delta(N_{\cal L}-\bigcup_jU(l_j))=
\Delta(N_{\cal L})\prod_{j=1}^n(t_j-1)$. Now, it is easily seen that 
$N_{\cal L}-\bigcup_jU(l_j)\cong N-U({\cal L})$, where $U({\cal L})$ denotes a 
tubular neighborhood the link. Moreover, the Alexander polynomial
of this is the ordinary one, 
$\Delta(N-U({\cal L}))=\Delta({\cal L})$, which yields the following
relation. 
\begin{equation}\label{eq-NLpoly}
\Delta(N_{\cal L})\;\;\;=\;\;\;\frac {\Delta({\cal L})}{(t_1-1)\ldots(t_n-1)}\;\; 
\end{equation}
We use (\ref{eq-NLpoly}) now to relate the $k=1$ torsion invariants $\Delta_{\varphi}$
to the Lescop invariant  $\lambda_L$ for 3-manifolds with $b_1(M)\geq 1$. 
Recall from \cite{Le96} that  $\lambda_L$ 
is an extension of the Casson-Walker invariant. 
\begin{lemma}\label{lm-Lescop}
 For a closed, compact, oriented manifold $M$ with $b_1(M)\geq 1$
and for \underline{any} epimorphism $\varphi:H_1(M)\onto \Z$ we have 
\begin{equation}\label{eq-Lescop}
\lambda_L(M)\;\;=\;\;\frac 12 \Delta''_{\varphi}(M)(1)-\frac 1 {12}
\Delta_{\varphi}(M)(1)\;. 
\end{equation}
\end{lemma}

\proof 
For  $b_1(M)=1$ Lescop gives in T5.1 in \S 1.5 of \cite{Le96} the formula
 $\lambda_L(M)=\frac 12 \Delta''(M)(1)-\frac 1 {12} |\Tors(M)|$. 
In this case $\varphi$ is canonical, and it follows from (\ref{eq-AlexForm})  that
$\Delta_{\varphi}(M)(1)=\Delta(M)(1)=|\Tors(M)|$.

For  $b_1(M)\geq 2$ define  the function 
$\zeta({\cal L})=\frac {\partial^n\Delta({\cal L})}
{\partial t_1\ldots\partial t_n}(1,\ldots,1)$ 
for a link ${\cal L}\subset N$ in a rational homology sphere 
with null homologous components and trivial linking matrix. In 2.1.2 of
\cite{Le96} Lescop then  writes $\lambda_L(N_{\cal L})=\zeta({\cal L})\,$.
 Combining this with (\ref{eq-NLpoly}) and the fact that every 3-manifold is
of the form $M=N_{\cal L}$, with $N$ and $\cal L$ as in (\ref{eq-NLpoly}), 
we thus have  $\lambda_L(M)=\Delta(M)(1,\ldots,1)$.  

Let now $\varphi:H_1(M)\onto \Z$ be an epimorphism for 
a closed 3-manifold, $M$, with $b_1(M)\geq 2$ 
so that we are  in the  situation (3) of 
Theorem~\ref{thm-turtor}. 
We know from our previous discussion that $\Delta_{\varphi}(M)(1)=0$ and 
that $\Delta_{\varphi}(M)$ is invariant under the substitution 
$t\symm t^{-1}$. This implies for the expansion in $(t-1)$ that 
$\Delta_{\varphi}(M)(t)=\frac 12 (t-1)^2\Delta''_{\varphi}(M)(1)\,+
\,{\cal O}((t-1)^3)\,$. By (3) of 
Theorem~\ref{thm-turtor} we therefore have 
$\overline{\varphi}(\Delta(M))(t)=\frac 1 2t  \Delta''_{\varphi}(M)(1)\,+
\,{\cal O}((t-1))\,$. Now, if $\varepsilon:\Z[\Pi]\onto\Z$ is the augmentation
map then $\lambda_L(M)=\Delta(M)(1,\ldots,1)=\varepsilon(\Delta(M))=
\varepsilon(\overline{\varphi}(\Delta(M)))=
\varepsilon(\frac 1 2t \Delta''_{\varphi}(M)(1)\,+\,{\cal O}((t-1)))
=\frac 12 \Delta''_{\varphi}(M)(1)\,$. 
\ep

The reason we find (\ref{eq-Lescop}) useful lies in the fact that
it does not distinguish between 
the cases $b_1(M)=1$ and $b_1(M)\geq 2$. Moreover,  any arbitrary 
$\varphi$ can be used to evaluate $\lambda_L$, which is, by construction, 
independent of the choice of $\varphi$. This will be essential when we
prove in Theorem~\ref{thm-lescut}, asserting that for a manifold $M$ with 
$\lambda_L(M)\neq 0$ we cannot remove more than one surface from $M$ 
without disconnecting $M$. 

It is now well known that for $b_1(M)\geq 1$  the Milnor torsion
$\tau(M)$ also equals the 
Seiberg-Witten invariant as shown in \cite{MT96}, \cite{Tu98}. In the case $b_1(M)=0$ 
 the Seiberg-Witten invariant is identified with a combination of the
Casson-Walker invariant and  $\tau(M)\in\Z[H_1(M)]$, see \cite{Ni01}. We will 
discuss the Casson invariant for integral homology spheres further
in Section~3 below.

\section{The Frohman-Nicas TQFT -- construction,  characters 
and hard-Lefschetz decompositions}

In \cite{FN92} Frohman and Nicas introduce a ($\Z/2\Z$-projective)
TQFT ${\cal V}_{\Z}^{FN}$
based on the intersection homology of 
the Jacobians $J(X)={\rm Hom}(\pi_1(X),U(1))$.
Specifically,
the functor ${\cal V}_{\Z}^{FN}$ associates to every surface $\Sigma$ the lattice
${\cal V}_{\Z}^{FN}(\Sigma)=H^*(J(\Sigma),\Z)\cong \ext * H_1(\Sigma,\Z)$. 
Furthermore, it assigns to a cobordism $C:\Sigma_A\to\Sigma_B$ 
a linear (lattice) map ${\cal V}_{\Z}^{FN}(C):  {\cal V}_{\Z}^{FN}(\Sigma_A)\to
{\cal V}_{\Z}^{FN}(\Sigma_B)$, which  
is computed (up to a sign)
from intersection numbers of the surface Jacobians with respect to a
Heegaard splitting of $C$. 
The  mapping class group acts canonically
on ${\cal V}_{\Z}^{FN}(\Sigma_g)$, factoring
through the symplectic quotient
$\Gamma_g\onto \Spl(H_1(\Sigma_g,\Z))=\Spl(2g,\Z)$.

As in \cite{FN92} we introduce an inner product $\langle\_,\_\rangle$
on $\ext *H_1(\Sigma_g)$ with respect to a complex structure   $J\in \Spl(2g,\Z)$\,, 
related to the symplectic skew form  $(\_,\_)$
by  $(x,Jy)=\langle x,y\rangle$ and $J^2=-1$. We fix  a standard 
homology basis  $a_1,\ldots, a_g,b_1,$ $\ldots,b_g$, which is an orthonormal as
well as a symplectic basis with $b_i=Ja_i$. Denote by
$\omega_g=\sum_ia_i\wedge b_i$ the standard  symplectic form.
In this setting ${\cal V}^{FN}_{\Z}$
assigns the map $\alpha\mapsto \alpha\wedge a_{g+1}$ 
to the cobordism that attaches the $(g+1)$-st 1-handle to a surface 
$\Sigma_g$. (Here the $a_{g+1}$-cycle is 
contractible into the interior of the cobordism.) The linear map for 
the cobordism that attaches 
the dual 2-handle is the respective conjugate. Hence, for the
standard handle body of genus $g$
in $S^3$ the assigned state is $\Omega_g=a_1\wedge\ldots \wedge a_g$,  
the volume form of the corresponding Lagrangian of the handle body. 
The complementary handle body is mapped to $\langle\Omega_g,\_\rangle$. 

In \cite{Ke00} we introduce the Lefschetz-$\sel_2$-action on $\ext *H_1(\Sigma_g)$
defined for the standard generators as $E.\alpha=\omega_g\wedge\alpha$,
$F=E^*$, and $\hat H.\alpha=({\rm deg}(\alpha)-g)\alpha$. It turns out that 
${\cal V}^{FN}_{\Z}$ is equivariant with respect to this action so that
we have the Lefschetz Decomposition, 
\begin{equation}\label{eq-FNLdec}
{\cal V}^{FN}_{\Z}\;\;\;=\;\;\;\bigoplus_{j\geq 1}V_j\otimes 
{\cal V}^{(j)}_{\Z}\;\;,
\end{equation}
where $V_j$ is the $j$-dimensional, simple representation of $\sel_2$. 
Each ${\cal V}^{(j)}_{\Z}$ 
is an irreducible, lattice TQFT, and 
can be defined  as the restriction  to 
$ker(F)\cap ker(\hat H+j-1)$. 
For a 3-manifold, $M$, and a given, $\varphi:H_1(M)\onto\Z$, let $\Sigma$ and
$C_{\Sigma}$ be a dual surface and a covering cobordism as in Section~1. 
As with any TQFT, it follows from $S$-equivalence, functoriality, and 
cyclicity of traces that the expressions
\begin{equation}\label{eq-momenta}
\Delta_{\varphi}^{(j)}(M)\;=\;\trace\Bigl({\cal V}^{(j)}_{\Z}(C_{\Sigma})\Bigr)
\;\;\in\Z\;
\end{equation}
do not depend on the choice $\Sigma$ but only on the pair $(M,\varphi)$. 
We will call the $\Delta^{(j)}_{\varphi}$ the {\em fundamental} torsion 
coefficients or fundamental
 torsion weights since they can be understood as characters of 
fundamental $\Spl(2g,\Z)$-representations. In  Theorem~4.4 of \cite{FN92}
Frohman and Nicas establish the relation of ${\cal V}^{FN}_{\Z}$ with  
the Alexander polynomial via a Lefschetz trace.
Combining their result with (\ref{eq-FNLdec}) and (\ref{eq-momenta})
we find the following expression:

\begin{theorem}[{\cite{FN92}, \cite{Ke00}}]\label{thm-LefTr}
Let $\varphi:H_1(M)\onto\Z$, $\Delta_{\varphi}$, $\Delta_{\varphi}^{(j)}$,
$\Sigma$, and  
$C_{\Sigma}:\Sigma\to\Sigma$ as above. 
\begin{equation}\label{eq-AlexPoly}
\Delta_{\varphi}(M)\;=\;
\trace\Bigl((-t)^{\hat H}{\cal V}^{FN}_{\Z}(C_{\Sigma})\Bigr)\;=\;\sum_{j\geq 1}
 [j]_{-t}\Delta^{(j)}_{\varphi}(M)\;. 
\end{equation}
\end{theorem}
We denote as usual $\,[n]_q=\frac {q^n-q^{-n}}{q-q^{-1}}\,\in\Z[q^{\pm 1}]$. 
For example,  (\ref{eq-AlexPoly}) implies that the fundamental coefficients 
$(\Delta^{(1)}, \Delta^{(2)},\ldots)$ for the $5_1$ and the $5_2$ knot are
$(-3,-2,0,\ldots)$ and $(0,1,1,0,\ldots)$ respectively.  
If we combine Theorem~\ref{thm-LefTr}  with Lemma~\ref{lm-Lescop} 
we find the following 
expression for the Lescop invariant in terms of the  
 $\Delta_{\varphi}^{(j)}$'s. 
\begin{theorem}\label{thm-LesMon} For a closed 3-manifold and any
 epimorphism $\varphi:H_1(M)\onto \Z$ the Lescop invariant is related 
to the fundamental torsion weights by 
\begin{equation}\label{eq-LesMon}
\qquad \qquad \lambda_L(M)\;\;=\;\;\sum_{j\geq 1} L^{(j)}\Delta^{(j)}_{\varphi}(M)\;\;,
\end{equation}
where 
$\displaystyle \,L^{(j)}
=(-1)^{j-1}\frac {j(2j^2-3)}{12}\,$. 
\end{theorem}
Clearly, 
the coefficients $L^{(j)}$ all lie in $\frac 1 {12}\Z$, and the first
few 
of them are given by 
$L^{(j)}=-\frac 1 {12}, -\frac 56, \frac {15}4, -\frac {29}3, 
 \frac {235}{12}, -\frac {69}2,\ldots\,$.
It is an interesting question whether there are further choices of
coefficients, other that the ones 
used in (\ref{eq-LesMon}),
 which would make the sum independent of the choice of $\varphi$. 
In \cite{Ke00} we also identify ${\cal V}_{\Z}^{FN}$
with the Hennings TQFT for the quasitriangular Hopf algebra 
$\Z/2\Z\ltimes\ext *\R^2$. This entails a calculus for determining 
$\Delta_{\varphi}(M)$, and hence also $\lambda_L$,
 from a surgery diagram, extending the traditional
Alexander-Conway Calculus for knots and links. It also allows us to remove the
sign-ambiguity of ${\cal V}_{\Z}^{FN}$ using the 2-framings on the cobordisms
that are the standard additional data in the quantum constructions.

\section{The Johnson-Morita homomorphisms and the\break Casson invariant}

The mapping class group representations
of the TQFT's ${\cal V}^{(j)}_{\Z}$ from the previous section all had
the Torelli groups $\,{\cal I}_g={\rm ker}(\Gamma_g\to \Spl(2g,\Z))\,$ in
their kernel. The goal of this section is to seek representations and
TQFT's with smaller kernels, and determine the information needed to 
describe the Casson invariant. 

Let us first recall the
construction of representations  of the mapping class groups 
$\Gamma_g$ from Section~7 of \cite{Ke01}, obtained from the ${\cal
V}^{(j)}_{\Z}$. They are indeed nontrivial on ${\cal I}_g$,
but they do vanish on
 the  subgroup ${\cal K}_g\subset {\cal I}_g$,
which is generated by Dehn twists along bounding cycles on the surface 
$\Sigma_g$.
 Denote also ${\cal K}_{g,1}$ and ${\cal I}_{g,1}$
the respective subgroups of $\Gamma_{g,1}$ for a surface with one boundary component.
 
In a series of papers 
Johnson and Morita have studied these subgroups and their
abelian quotients  extensively. For 
$H=H_1(\Sigma_g,\Z)$ let 
$\D U=\raise 2pt\hbox{$\ext 3 H$}\!\big /\!
\raise -2pt\hbox{$(\omega_g\wedge H)$}$. 
Johnson \cite{Jo80} constructs a homomorphism $\tau_1$ giving rise to the 
following $\Spl$-equivariant
short exact sequence. 
\begin{equation}\label{eq-Jon}
0\;\to\;{\cal K}_g\;\longrightarrow\;
{\cal I}_g\;
\stackrel{\mbox{$\tau_1$}}{
-\!\!\!-\!\!\!\longrightarrow}\; U\;\to\;0\;\;. 
\end{equation}
Here, $U$ is thought of as a free abelian group. In \cite{Mo93}  
Morita extends this to a homomorphism, $\tilde k\,$, 
on the entire mapping class group.
\begin{equation}\label{eq-JMhom}
0 \;\to\;{\cal K}_g\;\longrightarrow\;
\Gamma_g\;
\stackrel{\mbox{$\tilde k$}}{
-\!\!\!-\!\!\!\longrightarrow}
\;\frac 12 U\rtimes \Spl(2g,\Z)\;\to\;0\;\;.
\end{equation}
As in \cite{Ke01} we introduce for any $x\in\ext mH$ the maps
$\nu(x):\ext *H\to \ext {*+m}H$ and  $\mu(x):\ext *H\to \ext {*-m}H$,
given by $\nu(x).y=x\wedge y$ and $\mu(x)=\nu(Jx)^*$. It follows from
basic relations that $\mu(x)$ maps the $ker(F)$'s to each other,
decreasing the $\hat H$-weight by $m$, and that 
$\mu(\omega_g\wedge z)\bigl |_{ker(F)}=0$ for any $(m-2)$-form $z$. Thus, in 
the case $m=3$, $\mu$ factors for every $j\geq 1$ into an $\Spl$-equivariant map, 
\begin{equation}
\mu^{\flat}\;:\;\;U\;\longrightarrow\;
{\rm Hom}({\cal V}^{(j)}_{\Z}(\Sigma_g),{\cal V}^{(j+3)}_{\Z}(\Sigma_g))\;. 
\end{equation}
Any such map  serves as an extension map for  
a representation of $U\rtimes \Spl(2g,\Z)$ on
${\cal V}^{(j)}_{\Z}\oplus {\cal V}^{(j+3)}_{\Z}$, in which $U$ acts
non-trivially. We thus have extended $\Gamma_g$-modules
${\cal U}^{(j)}_{\Z}(\Sigma_g)$, which fit into a short exact sequence
as follows. 
\begin{equation}\label{eq-JMextZ}
0\;\to\;{\cal V}^{(j+3)}_{\Z}(\Sigma_g)\;\into \;
{\cal U}^{(j)}_{\Z}(\Sigma_g)\;\onto \; 
 {\cal V}^{(j)}_{\Z}(\Sigma_g) \;\to\;0\;\;\;. 
\end{equation}
The Casson invariant $\lambda_C$ for integral homology spheres is closely related 
to the subgroups ${\cal K}_g\subset {\cal I}_g$. Specifically,  
let us denote by $M_{\psi}$ the 3-manifold obtained via a 
Heegaard construction by cutting $S^3$ along a standard embedded
surface $\Sigma_g$ and repasting it with an element $\psi\in{\cal I}_g$. 
This yields the following assignment 
\begin{equation}\label{eq-Casmap}
\lambda^*\,:\;{\cal I}_g\,\longrightarrow\,\Z\;:\;\;\; 
\psi\,\mapsto\,\,\lambda^*(\psi):=\lambda_C(M_{\psi})\,.
\end{equation}
In \cite{Mo89}, \cite{Mo91} Morita studies this map thoroughly. Two important observations of his
are that  any $\Z$-homology sphere is of the form $M_{\psi}$
with $\psi\in{\cal K}_g\,$, 
and that $\lambda^*$ restricted to ${\cal K}_g$ is a homomorphism. 
Consequently, $\lambda^*$ is uniquely determined by its values on the generators
${\sf D}_{\cal C}\in{\cal K}_{g,1}$, given by Dehn twists along  bounding 
curves ${\cal C}\subset \Sigma_{g,1}$. The difficulty of this
description
of $\lambda^*$  lies in the fact that 
it is often not easy to find a presentation of a specific  homology sphere by a product
of ${\sf D}_{\cal C}$'s.

 The value of $\lambda^*$ on one of these generators of  ${\cal K}_g$
is obtained from the Alexander polynomial via 
$\lambda^*({\sf D}_{\cal C})=\frac 12 \Delta''_{\cal C}(1)$ considering $\cal C$
 as a 
knot in $S^3$. If $\Sigma_{h,\cal C}\subset \Sigma_{g,1}$ is a surface of genus $h$
bounded by $\cal C$, let 
$u_1,\ldots,u_h, v_1,\ldots,v_h$ be a symplectic basis of $H_1(\Sigma_{h,\cal C})$
and 
$\omega_{\cal C}=\sum_iu_i\wedge v_i$ the respective symplectic form of this surface. 
Moreover, for two homology cycles $a$
and $b$ in $\Sigma_g$ let $l_0(a,b)=lk(a,b^+)$ denote their linking number,
where $b^+$ denotes the ``push-off'' of the cycle $b$ in positive normal direction. 
Morita introduces a homomorphism 
\begin{equation}\label{eq-theta0}
\theta_0:\ext 2 H\otimes \ext 2 H: a\wedge b\otimes c\wedge d
\mapsto  l_0(a,c)l_0(b,d)-l_0(a,d)l_0(b,c)\end{equation}
and 
finds that $\lambda^*({\sf D}_{\cal C})=\theta_0(t_{\cal C})$, 
where $t_{\cal C}=-\omega_{\cal C}\otimes{\omega}_{\cal
C}$.

In the context  of TQFT interpretations it is now remarkable that $\theta_0$ 
can be reexpressed by matrix elements 
of operators on $\ext *H$ using the  multiplication and contraction maps $\nu$ and 
$\mu$ as before. More precisely, define a map
$\Psi:\ext * H\otimes \ext * H\longrightarrow {\rm End}(\ext *H)$ by
$\alpha\otimes \beta\mapsto \Psi(\alpha\otimes\beta)=\nu(\alpha)\circ\mu(\beta)\,$. 
Also, let $\Omega_g= a_1\wedge\ldots\wedge a_g\in\ext g H$ 
be the handle body state of the 
Frohman-Nicas theory as in Section~2. The next identity follows now from 
an exercise in multilinear algebra using the relations  in \cite{Ke01}. 
\begin{lemma}\label{lm-multlin} For any $\,A\in\,\ext 2 H\otimes \ext 2 H$ we have 
\begin{equation}\label{eq-multlin}
\theta_0(A)\;=\;\langle \Omega_g,\Psi(A)\Omega_g\rangle\;. 
\end{equation}
\end{lemma}
We introduce a restricted Lefschetz $\sel_2$-actions for 
the surface $\Sigma_{h,\cal C}$. Specifically, we have a subalgebra
 $\sel_2^{\,\cal C}$,  generated
by $E_{\cal C}=\nu(\omega_{\cal C})$, $F_{\cal C}=
J\circ E_{\cal C}^*\circ J^{-1} =\mu(\omega_{\cal C})$,
and $\hat H_{\cal C}=[E_{\cal C},F_{\cal C}]=-h+\sum_{i=1}^h
\nu(u_i)\nu(u_i)^*+\nu(v_i)\nu(v_i)^*\,$.
Moreover,  we introduce the 
standard quadratic Casimir Operator 
$Q_{\cal C}=E_{\cal C}F_{\cal C}+\frac 14 \hat H_{\cal C}(\hat H_{\cal C}-2)$
of $\sel_2^{\,\cal C}$, as well as $D_{\cal C}:= 
\frac 14 \hat H_{\cal C}(\hat H_{\cal C}-2)$. 
 
\begin{theorem}\label{thm-Cassmat} 
Let $\lambda^*$ be as in (\ref{eq-Casmap}),
$\Omega_g\in \ext g L$ the standard handle body state, 
${\sf D}_{\cal C}\in{\cal K}_g$ the Dehn twist along 
a bounding curve ${\cal C}$,
and $E_{\cal C}, F_{\cal C}, \hat H_{\cal C}, D_{\cal C}, Q_{\cal C}\in 
U(\sel_2^{\,\cal C})$ as above. Then 
\begin{equation}\label{eq-Cassmat}
\lambda^*({\sf D}_{\cal C})\;\;
=\;\;-\langle\Omega_g, E_{\cal C}F_{\cal C}.\Omega_g\rangle\;\;
=\;\;\langle\Omega_g, D_{\cal C}.\Omega_g\rangle
\,-\,\langle\Omega_g, Q_{\cal C}.\Omega_g\rangle
\;\;.
\end{equation}
Moreover, if $\psi\in\Gamma_g$ either preserves or reverses the Lagrangian 
decomposition $H=L\oplus L^{\perp}$ given by the standard handle body, then
\begin{equation}\label{eq-Cassymm}
\lambda^*({\sf D}_{\psi({\cal C})})\;=\;\lambda^*({\sf D}_{\cal C})\;. 
\end{equation} 
\end{theorem}
\proof The identity in (\ref{eq-Cassmat}) is readily obtained by 
combining Lemma~\ref{lm-multlin} and Morita's expression. 
The second assertion is obvious in the cases when $\psi$ preserves 
$L$ and $ L^{\perp}$, since 
the same is true for $\psi^{-1}$ and $\psi^*$ and we have 
$\bar\psi\circ X_{\cal C}\circ 
\bar\psi^{-1}=X_{\psi({\cal C})}$ for every generator 
$X\in\sel_2^{\,\cal C}$,  
where $\bar \psi$ the symplectic action on $\ext *H$.
 For the reversing case it now suffices to consider only a 
representative ${\sigma}\in \Gamma_g$ of the complex structure 
$\bar\sigma=J\in \Spl(2g,\Z)$.  
We observe that for all of the operators 
$Y_{\cal C}\in \{E_{\cal C}F_{\cal C}, Q_{\cal C}, D_{\cal C}\}\,$
we have  $Y_{\cal C}^*=JY_{\cal C}J^{-1}= Y_{\sigma(\cal C)}$, which 
implies (\ref{eq-Cassymm}).
\ep

The operator form of (\ref{eq-Cassmat}) suggests several 
decompositions of $\lambda^*$. 
Both $Q_{\cal C}$ and $D_{\cal C}$ have spectrum 
$\{\frac {j^2-1}4\,:\, j=1,2,\ldots,h'+1\}$, where 
$h'=\min(h,g-h)$.  To describe the eigenspaces more 
precisely, note that $H_1(\Sigma_g)\cong H_1(\Sigma_{h,\cal C})
\oplus  H_1(\Sigma_{g-h,\cal C})$, where $\Sigma_{g-h,\cal C}$ is
the complementary surface. Since $Q_{\cal C}$ also commutes with the
total $\sel_2$-action it preserves ${\cal V}^{(0)}_{\Z}(\Sigma_g)$,
which contains $\Omega_g$ and, hence, also $Q_{\cal C}.\Omega_g$. 
The $j$-th eigensubspace in this restriction is thus 
${\mathbb  E}_j(Q_{\cal C})=
{\cal V}^{(j)}(\Sigma_{h,\cal C})\otimes{\cal V}^{(j)}(\Sigma_{g-h,\cal C})$. 
Now, $D_{\cal C}$ does not commute with $\sel_2$ but still
preserves the total degree and hence maps $\ext gH_1(\Sigma_g)$ to itself. 
Correspondingly, we find in this restriction
${\mathbb  E}_j^{\pm}(D_{\cal C})=
\ext {h+1\pm j}H_1(\Sigma_{h,\cal C})\otimes
\ext {g-h-1\mp j}H_1(\Sigma_{g-h,\cal C})$. 

The Casson invariant is thus
expressible as a sum of terms 
$\frac {j^2-1}4\langle\Omega_g,{\mathbb P}^{(j)}_{\cal C}.\Omega_g\rangle $,
where the ${\mathbb P}^{(j)}_{\cal C}$ are 
projectors onto the eigenspaces of $Q_{\cal C}$ or $D_{\cal C}$.
Since $\bar\psi\circ {\mathbb P}^{(j)}_{\cal C}\circ 
\bar\psi^{-1}={\mathbb P}^{(j)}_{\psi({\cal C})}$  and 
$\Gamma_g$ acts transitively on the set of bounding curves,  
these operators and spaces can be determined from just the standard curve
${\cal C}_0^{(h)}$ for which $H_1(\Sigma_{{\cal C}_0^{(h)}})$ has
$a_1,\ldots,a_h,b_1,\ldots,b_h$ as a symplectic basis. 

One
well known decomposition of the Casson invariant is given by its
r\^ole in Floer cohomology. It is interesting to understand whether 
the eigenspace decompositions discussed here provide similar splittings
of  $\lambda^*$ or $\lambda_C$ as group morphisms or 3-manifold invariants
and are, possibly, even related to the Floer group decomposition.

Further investigations of Morita in \cite{Mo89}, \cite{Mo91} are devoted to understanding the 
failure of $\lambda^*$ to extend as a homomorphism to ${\cal I}_g$. More precisely,
he considers the integral cocycle (and rational coboundary) on  
${\cal I}_g$, given as
\begin{equation}\label{eq-Cascocy}
\delta_C(\phi,\psi)\;=\;\frac 12 \delta\lambda^*(\phi,\psi)\;=\;
\frac {\lambda^*(\phi\psi)-\lambda^*(\psi)-\lambda^*(\phi)} 2\;.
\end{equation} 
In Theorem 4.3 of \cite{Mo91} finds the expression 
$\delta_{\cal C}(\phi,\psi)=\tilde s(\tau_1(\phi),\tau_1(\psi))$, where $\tau_1$ is the
Johnson homomorphism onto $U=\ext 3 H/H$. The bilinear form $\tilde s(\ ,\ )$ descends
from the bilinear form $s$ on $\ext 3H$  defined by $s(\alpha,\beta)=(\alpha,\Pi_L\beta)$,
where $\Pi_L$ is the canonical projection onto $\ext 3 L\subset \ext 3H$ and $(\ ,\ )$
is 
the extension  of the standard symplectic form. It is easy to see that $\omega\wedge H$
lies in the left and right null space of $s$. From (\ref{eq-Cascocy}), 
Lemma~\ref{lm-multlin}, and further multilinear computations we obtain the following. 
\begin{theorem}\label{thm-cocycmat}
Let $\delta_{\cal C}$ be as in (\ref{eq-Cascocy}), 
$\tau_1$ as in (\ref{eq-Jon}), and $\mu$ and $\nu$ as before. Then 
\begin{equation}\label{eq-cocycmat}
\delta_C(\phi,\psi)\;=\;-\langle\Omega_g,\nu(\tau_1(\phi))\mu(\tau_1(\psi)).\Omega_g\rangle
\end{equation}
\end{theorem}  
Next, we describe a general form of a TQFT, motivated by the structure of
the Reshetikhin-Turaev Theories (see Section~4). This will provide a 
useful framework for finding TQFT interpretations of the  
 formulae  (\ref{eq-Cassmat}) and (\ref{eq-cocycmat}) for the 
Casson invariant.

We start with a commutative ring  $\M$ with unit, and denote the ring 
$\breve \M=\M[{\sf y}]/{\sf y}^2$ (of rank 2 as an $\M$-module). 
Moreover, we assume two $\Spl(2g,\Z)$-representations,  
$W_0(\Sigma_g)$ and $W_1(\Sigma_g)$, which are  free as modules over
$\M$. 
We write  $\brh_i:\Gamma_g\to {\rm GL}_{\M}(W_i(\Sigma_g))$ for the homomorphism 
with ${\cal I}_g$ in its kernel. Furthermore, we assume that each
$W_i$ admits an inner product $\langle\ ,\ \rangle$, and a special unit vector 
$\bOm_{g}\in W_0(\Sigma_g)$. We also denote 
 $\widetilde W_i(\Sigma_g)=W_i(\Sigma_g)\otimes_{\M} {\breve \M}$ and  
$\widetilde W(\Sigma_g)=\widetilde W_1(\Sigma_g)
\oplus \widetilde W_0(\Sigma_g)$ to which we extend $\langle\ ,\ \rangle$
with $W_1\perp W_0$. 
\begin{definition}\label{def-11sol} A {\em $1/1$-solvable TQFT} is a TQFT
 $\breve {\cal V}$ over a ring $\breve\M$ such that the 
$\breve\M$-modules are of the form 
$\breve {\cal V}(\Sigma_g)=\widetilde W(\Sigma_g)$ as above.
An element $\psi\in\Gamma_g$ is represented by $\breve {\cal V}$ in the form
\begin{equation}\label{eq-Vflp} 
\breve {\cal V}(\psi)\;=\;
\left[\begin{matrix} {\brh_1}(\psi) &\bmu (\psi)\cr 0 &
{\brh_0}(\psi)  \end{matrix}\right]
\;+\;{\sf y}\cdot \left[\begin{matrix} \bla_1(\psi)& \bka(\psi)\cr 
\bnu (\psi)& \bla_0(\psi) \end{matrix}\right]\;\;. 
\end{equation}
More generally, we require that the space 
$\widetilde W_1\oplus {\sf y}\cdot W_0$ 
is preserved (hence giving rise to a 
sub-TQFT over $\M$), 
and, furthermore, that $\breve {\cal V}$ assigns to the standard handle 
bodies the vectors
$\bOm$ and $\langle \bOm,\ \rangle$. 
Finally, the TQFT is {\em half-projective} with parameter 0 or ${\sf y}$. 
\end{definition}
Clearly, such a TQFT implies two invariants, $\tau^{\cal V}$
and $\lambda^{\cal V}$, of closed 3-manifolds into $\M$ defined as the polynomial
coefficients of the element in $\breve \M$ assigned by the TQFT 
as follows. 
\begin{equation}\label{eq-11solInv}
\breve {\cal V}(M)\;=\;\tau^{\cal V}(M)
\,+\,{\sf y}\cdot \lambda^{\cal V}(M)\;. 
\end{equation}
It also produces $\M$-valued invariants, $\Delta_{\varphi}^{\cal V}(M)$
and $\Xi_{\varphi}^{\cal V}(M)$,  of pairs $(M,\varphi)$, where 
$\varphi:H_1(M)\onto\Z$, $\Sigma$, and $C_{\Sigma}$ are as in Sections~1 and 2, 
by the following generalization of (\ref{eq-momenta}).
\begin{equation}\label{eq-11solTor}
\trace\Bigl(\breve{\cal V}(C_{\Sigma})\Bigr)\;=\;\Delta_{\varphi}^{\cal V}(M)
\,+\,{\sf y}\cdot \Xi_{\varphi}^{\cal V}(M)\;. 
\end{equation}
Next, let us record a number of immediate consequences of the above
definitions. 

\begin{lemma}\label{lm-22props} 
Let $\breve{\cal V}$ be a 1/1-solvable TQFT over a ring $\M$ with unit  for
which $2\in\M$ is not a zero-divisor. 
\begin{enumerate}
\item With the  boundary operator 
$\delta \bxi(\psi,\phi)=\brh_j(\psi)\bxi(\phi)- \bxi(\psi\phi)+
\bxi(\psi)\brh_i(\phi)$ for $\bxi:\Gamma_g\to {\rm Hom}(W_i,W_j)$ we have 
the relations 
\begin{eqnarray} 
\delta\bnu\;=\;\delta\bmu&=&0\label{eq-cohbnu}\\
-\delta\bla_1(\psi,\phi)&=&\bmu(\psi)\bnu(\phi)\label{eq-cohbla1}\\
 -\delta\bla_0(\psi,\phi)&=&\bnu(\psi)\bmu(\phi)\label{eq-cohbla0}\\
 -\delta\bka(\psi,\phi)&=&\bla_1(\psi)\bmu(\phi)+\bmu(\psi)\bla_0(\phi)\,
\label{eq-cohbka}
\end{eqnarray}
\item The restrictions of the maps $\bmu$ and $\bnu$ 
to ${\cal I}_g$ vanish on ${\cal I}_g'=[{\cal I}_g,{\cal I}_g]\,$. 
\item $\bmu$ and $\bnu$ also 
 factor through $\Spl(2g,\Z)$-equivariant, linear
maps $U\to {\rm Hom}_{\M}(W_0,W_1)$ and 
$U\to {\rm Hom}_{\M}(W_1,W_0)$ respectively.
\item The restrictions of the $\bla_i$ to  ${\cal K}_g$ vanish on
${\cal K}_g'=[{\cal K}_g,{\cal K}_g]$, and thus define  
$\Spl(2g,\Z)$-equivariant, linear
maps $H_1({\cal K}_g)\to {\rm End}_{\M}(W_i)$. 

\item For a Heegaard presentation $M_{\psi}$ we have 
\begin{equation}\label{eq-00mat}
\tau^{\cal V}(M_{\psi})=\langle\bOm,\brh_0(\psi)\bOm\rangle
\qquad 
\mbox{and}
\qquad
\lambda^{\cal V}(M_{\psi})=\langle\bOm,\bla_0(\psi)\bOm\rangle\;.
\end{equation}
\item If $M$ is a $\Z$-homology sphere then $\tau^{\cal V}(M)=1$.
\item The map 
$\lambda^{*\cal V}:\,\Gamma_g\to \M\,:\;\psi\mapsto
\lambda^{*\cal V}(\psi):=\lambda^{\cal V}(M_{\psi})$ restricts to 
a homomorphism on ${\cal K}_g$.
\item The cocycle $\delta\lambda^{*\cal V}(\psi, \phi)=-
\langle\bOm, \bnu(\psi)\bmu(\phi)\bOm\rangle\,$ restricted to ${\cal I}_g$
factors through a bilinear form on $U$. 
\end{enumerate}  
\end{lemma}

\proof The cohomological relations are an immediate consequence
of the fact that
the map $\breve{\cal V}:\Gamma_g\to {\rm GL}_{\breve\M}(\widetilde W(\Sigma_g))$ 
is a homomorphism. For example, (\ref{eq-cohbnu}) implies 
$\bmu({\cal I}_g')=0=\bnu({\cal I}_g')$.
>From Johnson's results, see Theorems~3 and 6 in
\cite{Jo85},  we have that ${\cal K}_g/{\cal I}_g'\,\cong (\Z/2\Z)^N$
so that also $\bmu({\cal K}_g)=0=\bnu({\cal K}_g)$ given that $(\M)_2=0$. 
Each of the remaining assertions follows now easily from previous assertions
and relations (\ref{eq-cohbnu})   and (\ref{eq-cohbka}). 
\ep

The functorial properties of $\lambda^{\cal V}$ are strikingly similar to
those of the Casson invariant $\lambda_C$. It is thus plausible to expect 
a TQFT-interpretation of $\lambda_C$ to come about in this form. 
We thus add the notion that a 1/1-solvable TQFT is of {\em Casson-type}
if $\overline{\lambda_C(M)}=\lambda^{\cal V}(M)$ for any $\Z$-homology sphere.
Here, we denote by $n\mapsto\overline n$ the canonical map $\Z\to\M$.  
The similarities of formulae is reflected in the following observation.
 
For $0<h<g$ and the  standard separating curve ${\cal C}_0^{(h)}$ as
before let us write $\bla^{(h)}:=\bla_0({\sf D}_{{\cal C}_0^{(h)}})$ and 
$L^{(h)}:=E_{{\cal C}_0^{(h)}}F_{{\cal C}_0^{(h)}}$. Note,  that both
operators act on $\Spl(2g,\Z)$-modules $W_0$ and $\ext *H$ respectively, and 
that they commute with the action of the standard subgroup
$\Spl(2h,\Z)\times \Spl(2(g-h),\Z)\subset \Spl(2g,\Z)$. The comparison 
of formula (\ref{eq-00mat}) in Lemma~\ref{lm-22props}  with 
(\ref{eq-Cassmat}) and (\ref{eq-cocycmat}) is summarized in the 
next lemma. 
\begin{lemma}\label{lm-22char}
A 1/1-solvable TQFT over $\Z$ is of Casson-type if and only if 
\begin{equation}\label{eq-prob1}
\langle \Omega_g, G\cdot L^{(h)}\cdot  G^{-1}.\Omega_g\rangle
\;=\; 
\langle \bOm_g, G\cdot  \bla^{(h)}\cdot  G^{-1}.\bOm_g\rangle
\end{equation}
for all $g,h\in\N$ with $0<h<g$, and for all
$G\in \raise 3pt 
\hbox{$\scriptstyle \Spl(2g,\Z)$}\!\big /\!
\raise -3pt 
\hbox{$\scriptstyle \Spl(2h,\Z)\times \Spl(2(g-h),\Z)$}$.
In this case we also have for all $\, a,b \in\,U\,$ the relation
 \begin{equation}\label{eq-prob2}
-2\langle\Omega_g, \nu(a)\mu(b).\Omega_g\rangle
\;=\; 
\langle\bOm_g, \bnu(a)\bmu(b).\bOm_g\rangle\;. 
\end{equation}
\end{lemma}
Finally, let us point out some 
subtleties associated with the second Johnson homomorphism. 
We write  $T\subset \ext 2 H\otimes \ext 2H$ for the symmetric subspace generated 
by $x\otimes x$ and $x\symm y:=x\otimes y+y\otimes x$ for all $x,y\in\ext 2H\,$, and,
further, denote by   
${\mathfrak h}_{g,1}(2)=T/T_0$ the quotient of $T$ 
by the subspace $T_0$ generated by elements
$a\wedge b\symm c\wedge d - a\wedge c\symm b\wedge d
+ a\wedge d\symm b\wedge c\,$, see \cite{Mo91}. 
The second Johnson homomorphism $\tau_2$ is now a map as follows. 
\begin{equation}\label{eq-Jon2}
\tau_2\;:\;\;{\cal K}_{g,1}\;\longrightarrow\;{\mathfrak h}_{g,1}(2)\;:\;\;\;
{\sf D}_{\cal C}\;\mapsto\;\overline{t_{\cal C}}\qquad
\mbox{with}\;\;\;t_{\cal C}:=-\omega_{\cal C}\otimes
\omega_{\cal C}\;.
\end{equation}
Here, $r\mapsto\overline r$ stands
for the map $T\to T/T_0$, and ${\sf D}_{\cal C}$ is
as in (\ref{eq-theta0}). Following (\ref{eq-theta0}), 
 the Casson invariant on ${\sf D}_{\cal C}$ is also given by the homomorphism
$\theta_0$ evaluated
on the element $t_{\cal C}\in T$. However, $\theta_0$ does not vanish on 
$T_0$ and thus
does not factor through ${\mathfrak h}_{g,1}(2)$. Consequently, 
$\lambda^*\big |_{{\cal K}_g}$ also does not factor through  $\tau_2$.  
 Yet,  in \cite{Mo91} Morita is able to define
a homomorphism $\eta:\Gamma_{g,1}\longrightarrow
{\mathbb Q}$ such that
$\eta({\sf D}_{\cal C})= \frac 16 h(h-1)\,$, as well as  a homomorphism
$\bar q_0:\,{\mathfrak h}_{g,1}(2)\longrightarrow
{\mathbb Q}$, such that $\lambda^*(\psi)=\eta(\psi)\,+
\,\bar q_0(\tau_2(\psi))\,$.  This raises the question what the 
relation is between  this decomposition and  the splitting of $\lambda^*$  
entailed by Theorem~\ref{thm-cocycmat}. Moreover, it is interesting to 
understand the r\^ole of $\tau_2$ in the general framework of  1/1-solvable
TQFT's.  

\section{$p$-Modular, Homological TQFT's -- their Relation to $S_n$ 
 resolutions, extensions and  characters} 

There are two ways to produce interesting TQFT's 
over the finite field $\F_p=\Z/p\Z$,
for a prime number $p\geq 3$. One is to consider 
the $\F_p$-reductions ${\cal V}^{(j)}_p$  of the Frohman-Nicas
lattice TQFT's ${\cal V}^{(j)}_{\Z}$. The second is obtained from the 
constant order reduction of the cyclotomic
integer expansion of the Reshetikhin-Turaev Theories. 
We will explore some relations between these two theories later.
As a preparation let us  first discuss
the properties of the ${\cal V}^{(j)}_p$'s 
and the $p$-modular versions of Theorems~\ref{thm-LefTr} and \ref{thm-LesMon}.

The ring reduction alone from $\Z$ to $\F_p$ turns the irreducible TQFT
${\cal V}^{(j)}_{\Z}$ into a generally highly reducible TQFT ${\cal V}^{(j)}_p$.
Specifically, the inner product on $\ext *H_1(\Sigma_g)$  induces  a pairing 
$\lz\_,\_\rz_p: {\cal V}^{(j)}_p(\Sigma_g)
\otimes {\cal V}^{(j)}_p(\Sigma_g) \to \F_p$. It is clear that the null-space
of this pairing yields a well-defined sub-TQFT.
\begin{definition}{\rm\cite{Ke01}}\label{def-quotTQFT}\qua
Let $\dov {\cal V}^{(j)}_p$ be the {quotient-TQFT} obtained from 
${\cal V}^{(j)}_p$ by dividing the vector space of each surface by
the null-space of $\lz\_,\_\rz_p$.
\end{definition}
Next, we illustrate explicitly that this is a nontrivial operation. 
\begin{example}
The map ${\cal V}_p^{(p-1)}(\Sigma_p)\onto \dov {\cal V}_p^{(p-1)}(\Sigma_p)$
has nontrivial kernel, given by $\F_p\overline{\omega_p}=im(E)$.
\end{example}
\proof Over $\Z$ the symplectic form $\omega_p$ is not in 
${\cal V}^{(p-1)}_{\Z}(\Sigma_p)$. However, we can pick another 
representative of $\overline{\omega_p}$, namely, 
$v=E.1-p(a_1\wedge b_1)=\omega_p - p(a_1\wedge b_1)\,\in\,\ext 2 H_1(\Sigma_p)$.
Since $F.a_i\wedge b_i=1$ we find $F.v=0$ so that indeed 
$v\in {\cal V}^{(p-1)}_{\Z}(\Sigma_p)=ker(F)\cap\ker(\hat H+p-2)$. 
Now, if $w\in{\cal V}^{(p-1)}_{\Z}(\Sigma_p)$ is any other such vector
we find $\lz v,w\rz=\lz E1,w\rz - \lz p(a_1\wedge b_1),w\rz =
\lz 1,Fw\rz - p\lz (a_1\wedge b_1),w\rz\,=\,- p\lz (a_1\wedge b_1),w\rz\,\in\;p\Z$.
Thus, if $\overline v =\overline {\omega_p}$ and $\overline w$ 
are the respective vectors in 
${\cal V}^{(p-1)}_{p}(\Sigma_p)$ we see that $\lz \overline {\omega_p},
 \overline w\rz_p=0$
so that $\overline{\omega_p}\neq 0$ lies in the null space of pairing $\lz\_,\_\rz_p$
and, hence, in the kernel of the above projection. The fact that the kernel
is not bigger than this is implied by Theorem~\ref{thm-res} below.
\ep

The general relation between the ${\cal V}_p$ and $\dov {\cal V}_p$ has the following
description. 
\begin{theorem}{\rm\cite{Ke01}}\label{thm-res}\qua The  TQFT's
$\dov {\cal V}^{(j)}_p$ are irreducible for any $j\in \N$ and any prime $p\geq 3$.
Each $\dov {\cal V}^{(j)}_p(\Sigma_p)$ carries a nondegenerate inner form, with
a compatible, irreducible $\Spl(2g,\Z)$-representation (i.e., $\psi^*=J\psi^{-1}J^{-1}$).

Moreover, for any $k\in \N$ with $0<k<p$ we have a resolution of the 
quotient-TQFT given by an exact sequence as follows.
\begin{equation}\label{eq-exactV}
\ldots\,\to\,{\cal V}^{(c_{i+1})}_p
\,\to\,{\cal V}^{(c_i)}_p\,\to\,\ldots
\,\to\, {\cal V}^{(2p+k)}_p
\,\to\,{\cal V}^{(2p-k)}_p
\,\to\,{\cal V}^{(k)}_p
\,\to\,\dov {\cal V}^{(k)}_p\,\to\,0\,, 
\end{equation}
where $c_i=ip+k$ if $i$ is even, and $c_i=(i+1)p-k$ if $i$ is odd.
\end{theorem}
The arrows in the sequence are to be understood as natural transformations
between TQFT functors. Particularly, this means that we obtain an
$\Spl(2g,\Z)$-equivariant resolution of $\dov {\cal V}^{(k)}_p(\Sigma_g)$ for
every $g\geq 0$, whose form (\ref{eq-exactV}) is, remarkably,  independent of $g$. 
 Quite curiously, the maps in (\ref{eq-exactV}) are given  by the Lefschetz 
operators from Section~2. More precisely, we prove in \cite{Ke01}
that for $j\equiv k\mod p$
the operator $E^k$ maps ${\cal V}^{(j)}_{\Z}(\Sigma_g)$ to 
${\cal V}^{(j-2k)}_{\Z}(\Sigma_g)\,+\,p\ext {g-j+2k+1}H_1(\Sigma_g)$.
Hence, we obtain well defined maps
$E^{k_i}:{\cal V}^{(c_{i})}_p\longrightarrow {\cal V}^{(c_{i-1})}_p$ in
the $\F_p$-reduction.

The rank reduction from  ${\cal V}_p$ to $\dov {\cal V}_p$ makes the representation
theory more challenging as well. The $\Spl(2g,\Z)$-representations for the integral
TQFT's ${\cal V}_{\Z}$ obviously lift to representations of $\Spl(2g,\R)$, and are,
therefore, highest weight representations in the sense of standard Lie theory. It 
follows from simple dimension counting that most of the $\dov {\cal V}_p$ cannot
be of such a form.

In \cite{Ke01} we prove exactness of (\ref{eq-exactV}) by breaking the sequence
down into the $\sp_{2g}$-weight spaces ${\cal W}^{(c)}_p(\varpi)$ for
a weight 
$\varpi$,  which are evidently preserved by the $E^{k_i}$-maps. 
Each ${\cal W}^{(c)}(\varpi)$ carries a natural, equivariant action
of the symmetric group $S_n$, where $n$ is the number of zero components 
of $\varpi$, see \cite{Ke01}. The ${\cal W}^{(c)}_p(\varpi)$ turn out to be 
isomorphic to the standard Specht modules ${\cal S}^{(c)}_{p}$ over $\F_p$
associated to the Young diagram $[\frac {n+c-1}2, \frac {n-c+1}2]$. The 
respective weight spaces of the $\dov {\cal V}^{(k)}_p$ are easily identified with
the irreducible $S_n$-modules ${\cal D}^{(k)}_p$ over $\F_p$, 
obtained, similarly, by an inner form reduction as in \cite{Ja78}. 
Exactness in (\ref{eq-exactV}),
is thus a consequence of the following result in the representation theory of
the symmetric groups.
\begin{theorem}{\rm\cite{Ke01}}\qua
Let ${\cal S}_p^{(c)}$ and ${\cal D}_p^{(k)}$ be $S_n$-modules as
above, and denote by $\chi^{(c)}$ and $\phi^{(k)}_p$ their characters, respectively.
We have a resolution as follows. 
\begin{equation}\label{eq-exactS}
\ldots\to {\cal S}^{(c_{i+1})}_p\,\to\,
{\cal S}^{(c_i)}_p\,\to\,\ldots\,\to\,
{\cal S}^{(2p-k)}_p\,\to\,{\cal S}^{(k)}_p\,\to\,
{\cal D}^{(k)}_p\,\to\,0
\end{equation}
Here, $k$ and the $c_i$ are as in Theorem~\ref{thm-res}. We obtain the 
relation 
\begin{equation}\label{eq-Schars}
\phi^{(k)}_p\,=\,\sum_{i\geq 0} (-1)^{i} \chi^{(c_i)}\;.
\end{equation}
\end{theorem}  
The proof uses the
 precise ordered modular structure of the ${\cal S}^{(c)}_p$ given by 
Kleshchev and Sheth in \cite{KS99}, which turns out to be sufficiently rigid to
prohibit any homology. The character expansion of irreducible $p$-modular
$S_n$-characters into $p$-reductions of the ordinary characters in
(\ref{eq-Schars}) is a direct consequence of  (\ref{eq-exactS}), and
appears to be  new in the modular representation theory of $S_n$. 

In order to extend the results from Section~2 to $\F_p$ we introduce, in 
analogy to (\ref{eq-momenta}), the
{\em $p$-modular}, fundamental torsion weights, given for 
 a pair $(M,\varphi:H_1(M)\onto \Z)$ by 
\begin{equation}\label{eq-pmodtor}
\dov\Delta^{(j)}_{\varphi,p}(M)\,=\,\trace(\dov {\cal
V}_{p}^{(j)}(C_{\Sigma}))\;\;\in
\;\F_p\;.
\end{equation}
The images of the Alexander polynomial and the Lescop
invariant 
in the  cyclotomic integers
are next expressed in the weights from (\ref{eq-pmodtor}). 

\begin{theorem}\label{thm-pTorcoef}
Let  $\fr_p:\Z[t,t^{-1}]\to\F_p[\zeta_p]$ be
the canonical  ring homomorphism, and denote, for $p\geq 5$, by 
$\overline L_{p}^{(j)}\in\F_p$ of the coefficients 
${L}^{(j)}\in \frac 1 {12}\Z$ from (\ref{eq-LesMon}). Then
\begin{eqnarray}\label{eq-pAlex} 
\fr_p\bigl(\Delta_{\varphi}(M)\bigr)\;&=&\;
\sum_{k=1}^{p-1}[k]_{-\zeta_p}\cdot\dov\Delta^{(k)}_{\varphi,p}(M)
\;\;\;\in\;\F_p[\zeta_p]\;.\\
\ \nonumber\\
\label{eq-LecMod}
\lambda_L(M)\;\;
&=&\;\;\sum_{k=1}^{p-1} \overline L^{(k)}_p\dov\Delta^{(k)}_{\varphi,p}(M)
\quad\mod p\;. 
\end{eqnarray}
\end{theorem}

\proof 
The resolution  from 
(\ref{eq-exactV}) implies, analogous to (\ref{eq-Schars}), 
the alternating series 
\begin{equation}\label{eq-modwexp} 
\dov\Delta^{(k)}_{\varphi,p}(M)=\sum_{i\geq 0}(-1)^i
\Delta^{(c_i)}_{\varphi}(M)
\,\mod p\qquad 0<k<p\,. 
\end{equation}
We note also that
$\fr_p([c_i]_{-t})=[c_i]_{-\zeta_p}=(-1)^i[k]_{-\zeta_p}=
(-1)^{i+k-1}[k]_{\zeta_p}$. Combining 
(\ref{eq-modwexp}) with  (\ref{eq-AlexPoly}) we 
obtain  the expansion (\ref{eq-pAlex})
of the Alexander polynomial in $\F_p[\zeta_p]$
in terms of  the irreducible, $p$-modular  weights from (\ref{eq-pmodtor}). 
For
 the $p$-reduction of the Lescop invariant in (\ref{eq-LecMod}) note that
for $p\geq 5$ we have $\frac 1 {12}\in\F_p$ so that the $p$-reductions 
$\overline L_p^{(j)}$ of the $L^{(j)}$ from (\ref{eq-LesMon}) 
are  well defined. For example,
$\overline L_{5}^{(j)}=2,0,0,2,0,\ldots$ and 
$\overline L_{7}^{(j)}=4,5,2,2,5,4,0,\ldots\,$. In general, we have 
$\overline L_{p}^{(p\pm j)}=\overline L_{p}^{(j)}$ so that
$\overline L^{(c_i)}_p=(-1)^i\overline L^{(k)}_p$. From this,
(\ref{eq-modwexp}),  and
(\ref{eq-LesMon}) we thus infer (\ref{eq-LecMod}). 
\ep

Finally,  let us note that the Johnson-Morita extension we constructed in 
(\ref{eq-JMextZ}) factors into the irreducible, $p$-modular quotients so
that we have for $0<k<p-3$ representations $\dov{\cal U}^{(k)}_{p}$ 
of   $U\rtimes \Spl(2g,\Z)$ over
$\F_p$, which represent $U$ nontrivially and admit  short exact sequences 
as follows.
\begin{equation}\label{eq-JMextp}
0\;\to\;\dov{\cal V}^{(k+3)}_{\Z}(\Sigma_g)\;\into \;
\dov{\cal U}^{(k)}_{p}(\Sigma_g)\;\onto \; 
\dov {\cal V}^{(k)}_{p}(\Sigma_g) \;\to\;0\;. 
\end{equation}

\section{The $\F_p[\zeta_p]$-Expansion of the Reshetikhin-Turaev\break TQFT,
 the structure of the Fibonacci case and cut-numbers}
 
Recall that the Reshetikhin-Turaev invariants 
for $U_{\zeta_p}(\so_3)$,
at a $p$-th root of unity $\zeta_p$, lie in the
cyclotomic integers $\Z[\zeta_p]$ if $p\geq 3$ is a prime.
Their expansions in ${\sf y}=(\zeta_p-1)$ 
yield the Ohtsuki-Habiro
invariants, which, in lowest order, 
are related to the previously discussed  torsion and Casson invariants. 
Gilmer \cite{Gi01} gives an abstract proof that 
the TQFT's ${\cal V}^{RT}_{\zeta_p}$ associated to $U_{\zeta_p}(\so_3)$ 
can be properly defined as TQFT's ${\cal V}^{I}_{\zeta_p}$ 
defined over the cyclotomic integers $\Z[\zeta_p]$ for a certain
restricted set of cobordisms. Consider the 
ring epimorphism $\Z[\zeta_p]\onto \F_p[{\sf y}]/{{\sf y}^{p-1}}:
\,{\zeta_p}\mapsto 1+{\sf y}$ as well as the 
TQFT ${\cal V}^{I}_{p,{\sf y}}$ over $\F_p[{\sf y}]/{{\sf y}^{p-1}}$ induced by it
from ${\cal V}^{I}_{\zeta_p}$. 
For given bases we can, therefore, consider
the  expansions in ${\sf y}$ of the linear map assigned to a cobordism 
$C$. They give rise to  further reduced TQFT's 
${\cal V}^{\supp j}_p$ over $\F_p[{\sf y}]/{{\sf y}^{j+1}}$ as follows. 
\begin{equation}\label{eq-expand}
{\cal V}^{I}_{p,{\sf y}}(C)\;=\;\sum_{k=0}^{p-2}
{\sf y}^k\cdot {\cal V}^{I}_{p,[k]}(C)
\quad\mbox{and}\quad{\cal V}^{\supp  j}_p\;:=\;\sum_{k=0}^{j}
{\sf y}^k\cdot {\cal V}^{I}_{p,[k]}(C)\;.
\end{equation}
Here, each ${\cal V}^{I}_{p,[k]}(C)$ is a matrix with entries in $\F_p$.  
We will focus below on the structure of the TQFT's
  ${\cal V}^{\supp 0}_p$  and ${\cal V}^{\supp 1}_p$ over $\F_p$ and  
$\breve \F_p=\F_p[{\sf y}]/{\sf y}^2$, respectively. 
\begin{conjecture}\label{con-p} Let $p\geq 5$ be a prime 
and $\displaystyle  q_p=\frac {p-3} 2$. 
\begin{enumerate}
\item[\rm A)] There are TQFT's, $\dov{\cal U}^{(k)}_p$, which extend the $\Gamma_g$
representations from (\ref{eq-JMextp}). 
\item[\rm B)] The TQFT ${\cal V}_p^{[\leq 0]}$ over ${\mathbb F}_p$ is a  
 quotient of sub-TQFT's of  the $q_p$-fold symmetric product 
$\mbox{\it \large S}^{\,q_p}\dov{\cal U}^{(1)}_p\,$. 
\item[\rm C)] The TQFT ${\cal V}_{\zeta_p}^I$ is {\em half-projective} with
parameter ${\sf x}=(\zeta_p-1)^{q_p}$, and, as such, has a 
``block-structured'' $\Z[\zeta_p]$-basis. 
\end{enumerate}
\end{conjecture}
The statements are not independent. Obviously B) only makes sense if A)
is true. Moreover, what we call a ``block-structure'' in C), meaning roughly
that integral bases can be obtained by sewing surfaces together, implies
the conditions  for 
half-projectivity. Recall from \cite{Ke98} that a half-projective
TQFT $\cal V$ over a ring ${\sf R}$
with respect to some ${\sf x}\in {\sf R}$ fulfills all of the 
usual TQFT axioms except for the following modification of the
functoriality with respect to compositions. 
For two cobordisms, $C_1$ and $C_2$, with well
defined composite $C_2\circ C_1$ we have 
${\cal V}(C_2\circ C_1)={\sf x}^{\mu(C_2,C_1)}
\cdot {\cal V}(C_2){\cal V}(C_1)$, where 
$\mu(C_2,C_1)\in\N\cup\{0\}$ is the rank of the connecting map
$H_1(C_2\cup C_1)\to H_0(C_2\cap C_1)$ in the respective Mayer-Vietoris
sequence. The tensor product axiom remains the usual. It is immediate
that  ${\cal V}$ satisfies the ordinary functoriality of TQFT's if 
we either restrict ourselves to connected surfaces or if ${\sf x}$ is 
invertible in ${\sf R}$  and we rescale ${\cal V}$. Non-semisimple 
Hennings TQFT's and homological 
gauge theories as in Section~2
are the first examples for ${\sf x}=0$ theories, see \cite{Ke00}. 
Following \cite{Ke98}, a consequence of Conjecture~\ref{con-p} C) is the following
conjecture raised first by Gilmer. 

\begin{conjecture}{\rm\cite{Gi01}}\label{con-Pat}\qua
For a closed, connected 3-manifold we have 
\begin{equation}
cut(M)\;\leq\; \frac 1 {q_p}{\mathfrak o}_p(M)\;.
\end{equation} 
\end{conjecture}
Recall that the  cut-number
$cut(M)$ of a closed, connected 3-manifold, $M$, is defined, alternatively, 
as the maximal number of components that a surface $\Sigma\subset M$
can have for which $M-\Sigma$ is still connected, or as the maximal rank of
a (non-abelian) free group $F$ such that there is an epimorphism 
$\pi_1(M)\onto F$. The quantum-order ${\mathfrak o}_p(M)$
of a closed 3-manifold is the maximal $k$ such that
${\cal V}^{RT}_{\zeta_p}(M)\in{\sf y}^k\Z[\zeta_p]\,$, where ${\sf y}=\zeta_p-1$ 
as before, assuming that we have a normalization for which  
${\cal V}^{RT}_{\zeta_p}(S^3)$ is a unit in $\Z[\zeta_p]$. 

The evidence for B) of Conjecture~\ref{con-p} comes, e.g, from explicitly matching dimensions for 
$g=1,2$ and general $p$, from comparison of the asymptotic behavior of the dimensions
as $g\to\infty$ (see \cite{Ke01}), and from further consistencies with the 
cyclotomic integer expansions. The form given in B) also implies that the
power $(I{\cal I}_g)^{q_p+1}\subset\Z[\Gamma_g]$ of the augmentation ideal of
the Torelli group is in the kernel of ${\cal V}_p^{[\leq 0]}$. In analogy to
Definition~\ref{def-11sol}, we will thus call ${\cal V}_p^{[\leq l]}$
a $q_p$/$l$-solvable theory. Murakami's result \cite{Mu95} together with B)
would hence imply that ${\cal V}_p^{[\leq 1]}$ is a $q_p$/1-solvable
{\em Casson}-TQFT in the generalized sense of Definition~\ref{def-11sol}.

The motivation and another strong piece of evidence for
Conjecture~\ref{con-p} is the following example.
\begin{theorem}{\rm\cite{KeF}}\label{thm-5}\qua
Conjecture~\ref{con-p} holds true for $p=5$. 
More precisely, we have the following isomorphism for the constant order TQFT. 
\begin{equation}\label{eq-FibId}
\dov {\cal U}^{(1)}_5\;\;\cong\;\;{\cal V}^{\supp 0}_5\;\;.
\end{equation}
\end{theorem}
Note first, that, with $\,q_5=1\,$, the product-TQFT 
in part B) of Conjecture~\ref{con-p} is simply 
$\dov{\cal U}^{(1)}_5$
itself, which is, clearly, consistent with (\ref{eq-FibId}). 

Moreover, the projective parameter of ${\cal V}^{RT}_{\zeta_5}$ simply 
becomes
${\sf x}={\sf y}=\zeta_5-1$. We sometimes call the Reshetikhin-Turaev TQFT
at a fifth root of unity the {\em Fibonacci  TQFT's} since the dimension
of ${\cal V}^{RT}_{\zeta_5}(\Sigma_g)$ is given, e.g., for even $g$ by 
$5^{\frac g2}f_{g-1}\,$,  where $f_0=0, f_1=1, f_2=1, \ldots$ are the  
Fibonacci numbers. (A similar formula holds for odd $g$, see \cite{Ke01}.) 
Note, that the Kauffman bracket skein theory associated to $U_{\zeta_p}(\so_3)$ has 
only two colors, namely 1 and $\rho$, subject to $\rho\otimes\rho=1\oplus\rho$.
Despite the seeming simplicity of the Fibonacci TQFT,  it is shown in \cite{FRW01, FRW00} 
to be of fundamentally greater complexity than the TQFT's for  $p=3,4,6$, which
are already interesting. The topological content of the $p=4$-TQFT, for 
example,  has been identified with the Rochlin 
invariant and the Birman-Craggs-Johnson homomorphisms \cite{Wr94}.

We state next the consequences of Theorem~\ref{thm-5} and
the identification in (\ref{eq-FibId}) that concern the explicit 
relations between the Casson-Walker-Lescop invariant $\lambda_{CWL}$
and the Fibonacci TQFT.

\begin{corollary}The TQFT \ ${\cal V}_5^{[\leq 1]}$
 over $\breve \F_5=\F_5[{\sf y}]/{\sf y}^2$ is a 1/1-solvable TQFT
over $\F_5$ of Casson-type
in the sense of Definition~\ref{def-11sol}, and, thus, defines the 
$\F_5$-valued invariants for closed 3-manifolds
${\tau}_5$ and $\lambda_5$ as in (\ref{eq-11solInv}), as well as the 
$\F_5$-invariants $\Delta_{\varphi,5}$ and $\Xi_{\varphi,5}$
for pairs $(M,\varphi)$ as in (\ref{eq-11solTor}). 
\begin{enumerate}
\item $\tau_5(M)=| H_1(M,\Z)|\,\mod\,5\,$ if $b_1(M)=0$ and 0 else wise. 
\item $\lambda_5(M)=\lambda_{CWL}(M)\,\mod\,5\,$. 
\item $\Delta_{\varphi, 5}(M)=-2\cdot \lambda_5(M)\,\mod\,5$. 
\item $\Xi_{\varphi,5}(M)\,=\,\Xi_5(M)$ is independent of $\varphi$.  
\end{enumerate}
\end{corollary}

\proof Recall that for closed manifolds ${\cal V}^{(j)}_{\Z}(M)=0$ if
$j\geq 2\,$, and ${\cal V}^{(1)}_{\Z}(M)={\cal V}^{FN}(M)$ yields the order of
the first homology group of $M$. The first claim thus follows from the
fact that ${\cal V}^{(1)}$ occurs precisely once in the resolution of 
${\cal V}_5^{[\leq 0]}$.

For the case $b_1(M)=0$ the identification with $\lambda_{CWL}$ follows
from Murakami's work \cite{Mu95}. In the Lescop case, $b_1(M)\geq 1$, we find from 
(\ref{eq-LecMod}) and (\ref{eq-FibId}) that 
$\lambda_L(M)=2(\dov\Delta^{(1)}_{\varphi,5}(M)+
\dov\Delta^{(4)}_{\varphi,5}(M))=2\cdot
(\trace(\dov {\cal V}^{(1)}_{5}(C_{\Sigma})
\oplus\dov {\cal V}^{(4)}_{5}(C_{\Sigma})))=
2\cdot \trace(\dov {\cal U}^{(1)}_{5}(C_{\Sigma}))=
2\cdot \trace({\cal V}^{[\leq 0]}_5(C_{\Sigma}))$. 
Now, it also follows from 
TQFT axioms that ${\cal V}_{\zeta_5}^I(M)=(\zeta_5-\zeta_5^{-1})\cdot 
\trace({\cal V}^{I}_{\zeta_5}(C_{\Sigma}))=2{\sf y}\cdot 
\trace({\cal V}^{[\leq 0]}_{5}(C_{\Sigma}))+{\cal O}({\sf y}^2)$. 
Comparison yields the assertion. The identity for 
$\Delta_{\varphi, 5}(M)$ can be read from this calculation as well. 
The claim for $\Xi_{\varphi, 5}(M)$ follows from the analogous identification
 with the next order 
Ohtsuki invariant that appears as the coefficient of ${\sf y}^2$.
\ep

Let us, moreover, comment on the consequences of Lemma~\ref{lm-22char} for
the Fibonacci case. In \cite{KeF} we show that $\bla^{(h)}$ is an orthogonal
 projector, whose kernel is naturally isomorphic to 
$\dov {\cal V}^{(1)}(\Sigma_h)\otimes \dov {\cal V}^{(1)}(\Sigma_{g-h})
\oplus \dov {\cal V}^{(4)}(\Sigma_h)\otimes \dov {\cal V}^{(4)}(\Sigma_{g-h})$.
For small genera this space can be related to the known eigenspaces of 
the summands of the operator 
$L^{(h)}\,=\,Q_{{\cal C}_0^{(h)}}-D_{{\cal C}_0^{(h)}}$, and, thereby, yields
an alternative proof of Murakami's result \cite{Mu95} in this rather special case.
Therefore, it seems likely
that with a better understanding of the structure for general
genera and primes it is possible to give an entirely independent proof of
the result in \cite{Mu95} based on purely representation theoretic and TQFT 
methods.

Finally, let us illustrate some concrete topological applications
of the structural theory presented in this article. Given Theorem~\ref{thm-5}
also  Conjecture~\ref{con-Pat}
becomes a theorem in the case $p=5$ as follows. 
\begin{theorem}{\rm\cite{GK}}\label{thm-cut5}\qua
\begin{equation}\label{eq-cut5}
cut(M)\;\leq\;{\mathfrak o}_5(M)\;. 
\end{equation}
\end{theorem}

Let us mention here two examples in which 
(\ref{eq-cut5}) allows us to determine the cut number $cut(M)$. 
Computations of this kind by classical means, generally, 
entail quite complicated and difficult problems in topology or group theory.

{\bf Example C1 } Consider the manifold $W$ obtained by 0-surgery
along the link with linking numbers 0 shown in the following figure. 

\begin{figure}[ht!]
\cl{\epsfxsize1.3in\epsfbox{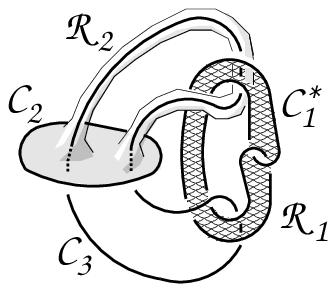}}
\end{figure}

We see that we have two disjoint surfaces that do not disconnect $W$.
They consist of the depicted Seifert surfaces ${\cal R}_i$ of ${\cal
C}_1^*$ and ${\cal C}_2$ and the discs glued in along the ${\cal
C}_i$'s by surgery.  Also, we know $b_1(W)=3$ so that $cut(W)$ can
still be either 2 or 3. It now follows from a short skein theoretic
calculation that ${\mathfrak o}_5(W)=2$ and hence $cut(W)=2$. 

\medskip
{\bf Example C2}\qua Let $\overline\psi\in
\Spl(H_1(\Sigma_g,\Z))$ be the symplectic, linear map associated to a
mapping class $\psi\in\Gamma_g$. 
We denote by $a_j(\psi)\in\Z$
the coefficients of the symmetrized, characteristic polynomial 
$t^{-g}\det(t{\mathbb I}_{2g}-\overline\psi)=\sum_ja_j(\psi)t^{j}$
so that $a_{-j}(\psi)=a_j(\psi)$ and $a_j(\psi)=0$ for $|j|>g$. 
Define next $\delta_5:\Gamma_g\to\F_5$ by 
$ \delta_5(\psi)=\sum_k a_{5k+2}(\psi)+a_{5k-2}(\psi)-a_{5k}(\psi)
\,\mod\,5$. Also, let $T_{\psi}=\Sigma\times[0,1]/\psi$ be the mapping torus 
for $\psi$. The combination of (\ref{eq-AlexForm}), (\ref{eq-pAlex}) from
Theorem~\ref{thm-pTorcoef}, 
(\ref{eq-FibId}), (\ref{eq-cut5}), and general TQFT
properties now yield the following criterion. 
\begin{lemma}If $\delta_5(\psi)\neq 0$ then 
$\,cut(T_{\psi})=1\,$.
\end{lemma}
Note here that the left hand condition only depends on the action
$\overline\psi$ on homology, and, e.g.,  in the case  $g=2$ 
reduces to 
$\trace(\overline\psi^2)+1\not\equiv \trace(\overline\psi)^2\mod 5\,$. 
(For $g=1$ we always have $cut(T_{\psi})=1$.)
 The precise knowledge of the higher order structure
of ${\cal V}_{\zeta_5}^I$ allows for finer theorems of this type,
and we, obviously, expect similar results to hold for general $p$.
\medskip

The last example is in fact a special case of a more general relation
between the Lescop invariant and cut-numbers, which is independent 
of the Reshetikhin-Turaev Theory. 

\begin{theorem}\label{thm-lescut} 
Let $M$ be a 3-manifold as before with $b_1(M)\geq 1$. Then, 
if $cut(M)\geq 2\,$, $\lambda_L(M)=0\,$. 
\end{theorem}

\proof The condition $cut(M)\geq 2$ means that $M-\Sigma-\Sigma'$ is 
connected for two embedded, oriented, two-sided surfaces, which means that
$C_\Sigma-\Sigma'$ is connected. This implies that $C_\Sigma-\Sigma'-\Sigma''$
splits into exactly two connected components for a surface $\Sigma''\neq\emptyset$.
Thus $C_\Sigma=A\circ B$ with connected cobordisms 
$A:\Sigma'\sqcup\Sigma''\to\Sigma$ and 
$B:\Sigma\to\Sigma'\sqcup\Sigma''\,$ so that $\mu(A,B)\geq 1$. As a result
of ${\sf x}=0$ half-projectivity of the Frohman-Nicas TQFT we thus
obtain ${\cal V}^{FN}(C_{\Sigma})=0$ and hence $\Delta_{\varphi}(M)=0$
for $\varphi$ dual to $\Sigma$. By (\ref{eq-Lescop}) this now implies
$\lambda_L(M)=0$. 
\ep 

In \cite{GK} we will give examples of $T_{\psi}$, with $\psi\in{\cal I}_g$, for
which $cut(T_{\psi})$
can no longer be determined from $\lambda_L$ or Alexander polynomials, but
where we have to employ Theorem~\ref{thm-cut5} to determine cut-numbers greater or
equal to 2.

\Addresses\recd

\end{document}

%% file: gtmonout.tex
\def\ifplaintex{\expandafter\ifx\csname documentclass\endcsname\relax}

\ifplaintex 
\else
\fi

\def\gtm{{\mathsurround=0pt\it $\cal G\mskip-2mu$eometry \&\ 
$\cal T\!\!$opology $\cal M\mskip-1mu$onographs}}    

\def\gtp{{\mathsurround=0pt\it $\cal G\mskip-2mu$eometry \&\ 
$\cal T\!\!$opology $\cal P\!$ublications}}  

\def\recd{{\small Received:\qua\receiveddate\ifx\reviseddate\relax
\else\qquad Revised:\qua\reviseddate\fi\par}} 

\def\volumenumber#1{\def\thevolumenumber{#1}}
\def\volumeyear#1{\def\thevolumeyear{#1}}
\def\volumename#1{\def\thevolumename{#1}}
\def\papernumber#1{\def\thepapernumber{#1}}
\def\pagenumbers#1#2{\def\startpage{#1}\def\finishpage{#2}}
\def\published#1{\def\publishdate{#1}}
\def\received#1{\def\receiveddate{#1}}
\def\revised#1{\def\reviseddate{#1}}
\def\accepted#1{\def\accepteddate{#1}}
\def\asciititle#1{\def\theasciititle{#1}}

\def\asciikeywords#1{\def\theasciikeywords{#1}}

\let\\\par
\let\thevolumenumber\relax\let\thepapernumber\relax
\let\thevolumeyear\relax\let\startpage\relax
\let\finishpage\relax\let\publishdate\relax\let\receiveddate\relax
\let\reviseddate\relax\let\accepteddate\relax\let\theasciititle\relax
\let\theasciiauthors\relax
\let\theasciiabstract\relax\let\theasciikeywords\relax

\let\theerratum\relax\let\theasciiemail\relax
\let\theshortauthors\relax\let\theshorttitle\relax

\def\startpage{1}\def\finishpage{15}\def\thepapernumber{77}

\volumenumber{2}
\volumename{Proceedings of the Kirbyfest}
\volumeyear{1999}

\long\def\maketitlep{   

\count0=\startpage

\gtm\nl        
{\small Volume \thevolumenumber: \thevolumename\nl 
\ifx\theerratum\relax\else Erratum \erratumnumber\nl\fi
Pages \startpage--\finishpage\nl}

\vglue 0.1truein   

{\parskip=0pt\leftskip 0pt plus 1fil\def\\{\par\smallskip}{\ifplaintex\large
\else\Large\fi\bf\thetitle}\par\medskip}   
\vglue 0.05truein 

{\parskip=0pt\leftskip 0pt plus 1fil\def\\{\par}{\sc\theauthors}
\par\medskip}%
 
\vglue 0.03truein 

{\small\leftskip 25pt\rightskip 25pt{\bf Abstract}\stdspace\theabstract

{\bf AMS Classification}\stdspace\theprimaryclass
\ifx\thesecondaryclass\relax\else; \thesecondaryclass\fi\par
{\bf Keywords}\stdspace \thekeywords\par}\vglue 7pt

}   

\font\phead=cmsl9 scaled 950
\font=cmsl9 scaled 1050
\font\pnum=cmbx10 scaled 913
\font\lnum=cmbx10 
\font\pfoot=cmsl9 scaled 950
\font=cmsl9 scaled 1050
\ifplaintex
\headline{\vbox to 0pt{\vskip -4.5mm\line{\small\phead\ifnum
\count0=\startpage ISSN 1464-8997 (on line)
1464-8989 (printed) \hfill {\pnum\folio}\else\ifodd\count0\def\\{ }%
\ifx\theshorttitle\relax\thetitle\else\theshorttitle\fi\hfill{\pnum\folio}
\else\def\\{ and }{\pnum\folio}\hfill\ifx\theshortauthors\relax\theauthors
\else\theshortauthors\fi\fi\fi}\vss}}
\footline{\vbox to 0pt{\vglue 0mm\line{\small\pfoot\ifnum\count0=\startpage
Published \publishdate:\qua\copyright\ \gtp\hfill\else
\gtm, Volume \thevolumenumber\ (\thevolumeyear)\hfill\fi}\vss
}}
\else
\makeatletter
\def\@oddhead{{\small\lhead\ifnum\count0=\startpage ISSN 1464-8997 (on line)
1464-8989 (printed) \hfill {\lnum\number\count0}\else\ifodd\count0
\def\\{ }\ifx\theshorttitle\relax \thetitle \else\theshorttitle\fi\hfill
{\lnum\number\count0}\else\def\\{ and }{\lnum\number\count0}
\hfill\ifx\theshortauthors\relax 
\theauthors\else\theshortauthors\fi\fi\fi}}\def\@evenhead{@oddhead}
\def\@oddfoot{\small\lfoot\ifnum\count0=\startpage Published \publishdate:\qua\copyright\ \gtp\hfill\else
\gtm, Volume \thevolumenumber\ (\thevolumeyear)\hfill\fi}
\def\@evenfoot{@oddfoot}
\makeatother
\fi

\let\maketitlepage\maketitlep
\let\makeshorttitle\maketitlepage
\let\maketitle\maketitlepage

\newwrite\gtoutfile
\long\gdef\makeheadfile{  
{\def\\{, }\def\s{ }
\immediate\openout\gtoutfile head.xxx
\immediate\write\gtoutfile{To: math@arxiv.org}
\immediate\write\gtoutfile{Subject: put OR rep NNNNN:ppppp}
\immediate\write\gtoutfile{--text follows this line--}
\immediate\write\gtoutfile{Proxy-for: \ifx\theasciiauthors\relax
\theauthors\else\theasciiauthors\fi\s<\ifx\theasciiemail\relax\theemail\else\theasciiemail\fi>}
\immediate\write\gtoutfile{\noexpand\\}
\immediate\write\gtoutfile{Authors: \ifx\theasciiauthors\relax
\theauthors\else\theasciiauthors\fi}
{\def\\{ }\immediate\write\gtoutfile{Title: \ifx\theasciititle\relax
\thetitle\else\theasciititle\fi}}
\immediate\write\gtoutfile{Subj-class: GT or SG, GR etc}
\immediate\write\gtoutfile{MSC-class: \theprimaryclass\ifx\thesecondaryclass\relax\else, \thesecondaryclass\fi}
\immediate\write\gtoutfile{Journal-ref: Geom. Topol. Monogr. \thevolumenumber\s
(\thevolumeyear) \startpage-\finishpage}
\immediate\write\gtoutfile{Comments: Published by Geometry and Topology Monographs at}
\immediate\write\gtoutfile{\s\s\s  http://www.maths.warwick.ac.uk/gt/GTMon\thevolumenumber/paper\thepapernumber.abs.html}
\immediate\write\gtoutfile{\noexpand\\}
\immediate\write\gtoutfile{}
\ifx\theasciiabstract\relax
\immediate\write\gtoutfile{\theabstract}\else
\immediate\write\gtoutfile{\theasciiabstract}\fi
\immediate\write\gtoutfile{}
\immediate\write\gtoutfile{\noexpand\\}
\immediate\write\gtoutfile{}
\immediate\closeout\gtoutfile}}  

\def\maketitlepage{\maketitlep\makeheadfile}
\let\makeshorttitle\maketitlepage
\let\maketitle\maketitlepage